\documentclass[12pt]{amsart}
\usepackage{amsmath}
\usepackage{amssymb}

\newtheorem{Thm}{Theorem}[section]
\newtheorem{Lemma}[Thm]{Lemma}
\newtheorem{Prop}[Thm]{Proposition}
\newtheorem{Def}[Thm]{Definition}
\newtheorem{Cor}[Thm]{Corollary}
\newtheorem{Rem}[Thm]{Remark}

\newtheorem{Example}[Thm]{Example}

\newtheorem{Question}{Question}

\newcommand{\e}{\varepsilon}

\newcommand{\dist}{\textrm{dist}}

\newcommand{\Lc}{L}
\newcommand{\dns}{\textrm{DN-S}}
\newcommand{\dss}{\textrm{DSS}}

\newcommand{\N}{\mathbb{N}}

\begin{document}

\title[Disjointly non-singular operators]{Disjointly non-singular operators on Banach lattices}

\thanks{Supported in part by MICINN (Spain), Grant PID2019-103961GB-C22.
}

\author[M.\ Gonz\'alez]{Manuel Gonz\'alez}
\address{Departamento de Matem\'aticas, Facultad de Ciencias,
Universidad de Canta\-bria, E-39071 Santander, Spain}
\email{manuel.gonzalez@unican.es}

\author[A.\ Mart\'\i nez-Abej\'on]{Antonio Mart\'\i nez-Abej\'on}
\address{Departamento de Matem\'aticas, Facultad de Ciencias,
Universidad de Oviedo, E-33007 Oviedo, Espa\~na}
\email{ama@uniovi.es}

\author[A.\ Martin\'on]{Antonio Martin\'on}
\address{Departamento de An\'alisis Matem\'atico, Facultad de Ciencias,
Universidad de La Laguna, E-38271 La Laguna (Tenerife), Spain} \email{anmarce@ull.es }


\begin{abstract}
An operator $T$ from a Banach lattice $E$ into a Banach space is disjointly non-singular
($\dns$, for short) if no restriction of $T$ to a subspace generated by a disjoint sequence
is strictly singular.
We obtain several results for $\dns$ operators, including a perturbative characterization.
For $E=L_p$ ($1< p<\infty$) we improve the results, and we show that the $\dns$ operators
have a different behavior in the cases $p=2$ and $p\neq 2$.
As an application we prove that the strongly embedded subspaces of $L_p$ form an open
subset in the set of all closed subspaces.
\end{abstract}

\maketitle

\thispagestyle{empty}

\section{Introduction}

Inspired by the study of tauberian operators $T:L_1\to Y$ using Banach lattice techniques
in \cite{GM:97}, we introduce the \emph{disjointly non-singular operators} from a Banach
lattice $E$ into a Banach space $Y$, denoted $\dns(E,Y)$, and the \emph{dispersed subspaces}
of $E$.
Note that $T:L_1\to Y$ is tauberian if and only if $T\in \dns(L_1,Y)$.

We give several characterizations of the operators in $\dns(E,Y)$ in terms of the action on
the normalized disjoint sequences in $E$.
We show that an operator $K:E\to Y$ is disjointly strictly singular \cite{Hdez-Salinas:89}
if and only if $M$ is dispersed when the restriction $K|_M$ is an isomorphism, that each
$T\in\dns(E,Y)$ has dispersed kernel $N(T)$, and that $T$ in $\dns$ and $K$ disjointly
strictly singular imply $T+K$ in $\dns$.
We prove a \emph{perturbative characterization:} $T:E\to Y$ is in $\dns$ if and only if
$N(T+K)$ is dispersed for every compact operator $K:E\to Y$, and we show that $T:E\to Y$ with
closed range is in $\dns$ if and only if $N(T)$ is dispersed.

To study the operators in $\dns(L_p,Y)$ ($1< p<\infty$), we show that the dispersed subspaces
of $L_p$ ($1\leq p<\infty$) coincide with the strongly embedded subspaces
\cite[Definition 6.4.4]{A-Kalton:06}, and that these subspaces form an open subset in the set
of all closed subspaces of $L_p$ with respect to the gap metric.
Since for $p\neq 2$ a closed subspace of $L_p$ is strongly embedded if and only if it does not
contain copies of $\ell_p$ \cite[Theorems 6.4.8 and 7.2.6]{A-Kalton:06}, but this is not true
for $p=2$, these two cases have a different behavior.
In fact, there are decompositions of $L_2$ as a direct sum of two strongly embedded subspaces
\cite[Theorem 8.22]{Pisier:CBMS}, but this is not true for $p\neq 2$.
We show that these decompositions of $L_2$ are stable under small perturbations with respect to
the gap metric.
Moreover $\dns(L_p)$ is stable under products for $p\neq 2$, but $\dns(L_2)$ is not.
We also prove that $\dns(L_p,Y)$ is non-empty if and only if $Y$ contains a subspace isomorphic
to $L_p$, that the disjointly strictly singular operators form the perturbation class of
$\dns(L_p,Y)$, and that there exist non-trivial examples of operators in $\dns(L_p)$ for
$1<p<\infty$.
\medskip

Throughout the paper, $E$ and $F$ will be Banach lattices and $X$, $Y$ and $Z$ will be Banach
spaces.
We denote by $\Lc(X,Y)$ the space of all (continuous linear) operators from $X$ into $Y$ and we
write $\Lc(X)$ when $X=Y$, $L_p$ is the space $L_p(0,1)$ for $1\leq p\leq \infty$, and we denote
by $[x_n]\equiv [x_n :n\in\N]$ the closed subspace generated by a sequence $(x_n)$ in $X$.

An operator $T\in\Lc(X,Y)$ is \emph{upper semi-Fredholm} if it has closed range $R(T)$ and finite
dimensional kernel $N(T)$.
The operator $T$ is \emph{strictly singular} if there is no closed infinite dimensional subspace
$M$ of $X$ such that the restriction $T|_M$ is an isomorphism.
The operator $T$ is \emph{$Z$-singular} if there exists no subspace $M$ of $X$ isomorphic to $Z$
such that the restriction $T|_M$ is an isomorphism.
The $\ell_1$-singular operators form an operator ideal and they are called \emph{Rosenthal operators}
by Pietsch in \cite{Pietsch}.

An operator $T\in\Lc(X,Y)$ is called \emph{tauberian} if the second conjugate
$T^{**}\in\Lc(X^{**}, Y^{**})$ satisfies $(T^{**})^{-1}(Y) = X$.
Introduced by Kalton and Wilansky in \cite{KaltonWil:76}, these operators are useful in several
topics of Banach space theory like factorization of operators, equivalence between the KMP and the
RNP, and refinements of James' characterization of reflexive Banach spaces.
We refer to \cite{GM:10} for information on tauberian operators.

\section{Disjointly non-singular operators}

We say that a sequence $(x_n)$ in $E$ is \emph{almost disjoint} if there exists a normalized
disjoint sequence $(y_n)$ in $E$ such that $\lim_{n\to\infty}\|x_n-y_n\|=0$.

\begin{Def}
A closed subspace $M$ of $E$ is said to be \emph{dispersed} if it does not contain almost
disjoint sequences.
\end{Def}

Recall that each disjoint sequence of non-zero vectors in $E$ is an unconditional basic sequence.

\begin{Lemma}\label{lemma:fin}
Let $M$ be a closed subspace of $X$ and let $(x_n)$ be a basic sequence in $X$.
\begin{enumerate}
  \item If $M\cap [x_n]$ is finite dimensional, then $M\cap [x_n : n> n_0]=\{0\}$ for some
  $n_0\in\N$.
  \item If, additionally, $M+ [x_n]$ is closed then so is $M+ [x_n : n> n_0]$.
\end{enumerate}
\end{Lemma}
\begin{proof}
(1) Let $(x^*_n)$ be a sequence in $X^*$ such that $x^*_i(x_j)=\delta_{i,j}$ for $i,j\in\N$
and suppose that $\dim M\cap [x_n]=k$.
If $0\neq x\in M\cap [x_n]$ there exists $l\in\N$ such that $x^*_l(x)\neq 0$.
Thus $x\notin M\cap [x_n :n>l]$, hence $\dim M\cap [x_n :n>l]<k$, and repeating the process
we get the desired $n_0\in\N$.

(2) Note that $M+[x_n:n>n_0]$ is the sum of $(M+[x_n])\cap \left(\cap_{k=1}^{n_0}N(x^*_k)\right)$
and a finite dimensional subspace.
\end{proof}

The argument in the proof of the following result will be used several times.

\begin{Lemma}[Perturbation lemma]\label{lemma:perturba}
Let $M$ and $N$ be closed subspaces of $X$, and let $0<\e<1/2$.
If $M+N$ is not closed, then there exists a compact operator $K\in\Lc(X)$ with $\|K\|<\e$
and a closed infinite dimensional subspace $M_1$ of $M$ such that $(I-K)(M_1)\subset N$,
$I-K$ is invertible and $(I-K)^{-1}=I+L$ with $\|L\|<2\e$.
\end{Lemma}
\begin{proof}
Since the sum $M+N$ is not closed, both subspaces $M$ and $N$ are infinite dimensional, and we can
select a normalized sequence $(m_k)$ in $M$, and sequences $(n_k)$ in $N$ and $(x^*_k)$ in $X^*$ so
that $x^*_i(m_j)=\delta_{i,j}$ and $\|x^*_k\|\, \|m_k-n_k\|<\e 2^{-k}$ for $i,j,k\in\N$.
See the proof of \cite[Proposition]{GO:86}.

Thus $Kx=\sum_{k=1}^\infty x^*_k(x) (m_k-n_k)$ defines a compact operator $K\in\Lc(X)$ such that
$(I-K)m_i=n_i$ for each $i\in\N$ and $\|K\| \leq \sum_{k=1}^\infty \|x^*_k\|\, \|m_k-n_k\|<\e$.
Hence $I-K$ is an invertible operator and we can take $M_1=[m_k]$.
\medskip

For the last part, apply $(I-K)^{-1} = I +\sum_{n=1}^\infty K^n$.
\end{proof}

Next we characterize the dispersed subspaces.

\begin{Prop}\label{prop:dispersed}
A closed subspace $M$ of $E$ is dispersed if and only if, for every disjoint sequence of
non-zero vectors $(x_n)$ in $E$, $M\cap [x_n]$ is finite dimensional and $M + [x_n]$ is closed.
\end{Prop}
\begin{proof}
Suppose that there is a disjoint sequence of non-zero vectors $(x_n)$ in $E$ such that
$M\cap [x_n]$ is infinite dimensional, or $M\cap [x_n]$ finite dimensional and $M + [x_n]$ is
not closed.
In both cases we can construct a normalized block-basis $(y_k)$ of $(x_n)$ with
$\lim_{k\to\infty}\dist(y_k,M)=0$.
Since $(y_k)$ is a disjoint sequence, $M$ is not dispersed.

Conversely, suppose that for every disjoint sequence of non-zero vectors $(x_n)$ in $E$,
$M\cap [x_n]$ is finite dimensional and $M + [x_n]$ is closed.
Given a normalized disjoint sequence $(x_n)$ in $E$, by Lemma \ref{lemma:fin} we have
$M\cap [x_n : n\geq n_0]=\{0\}$ for some $n_0$.
Then $\liminf_{n\to\infty}\dist(x_n,M)>0$; hence $M$ contains no almost disjoint sequence.
\end{proof}

The following class of operators was introduced in \cite{Hdez-Salinas:89}.

\begin{Def}
An operator $T\in\Lc(E,Y)$ is said to be \emph{disjointly strictly singular} if there is no
disjoint sequence of non-zero vectors $(x_n)$ in $E$ such that the restriction of $T|_{[x_n]}$
is an isomorphism.
\end{Def}

We denote $\dss(E,Y):=\{T\in\Lc(E,Y) : T \textrm{ is disjointly strictly singular}\}$.
\medskip

For $1\leq q <p<\infty$, the natural inclusion $L_p\to L_q$ is in $\dss$
\cite{Hdez-Salinas:89}, but it is not strictly singular because $\|\cdot\|_p$ and
$\|\cdot\|_q$ are equivalent in the subspace generated by the Rademacher functions.
Also, $\dss(E,Y)$ is a closed subspace of $\Lc(E,Y)$ \cite{Hdez:90}.

\begin{Prop}\label{prop:DSS-SS}
For $T\in\Lc(E,Y)$, the following assertions are equivalent:
\begin{enumerate}
  \item $T$ is disjointly strictly singular.
  \item Given a closed subspace $M$ of $E$, if $T|_M$ is an isomorphism then $M$ is dispersed.
  \item For every disjoint sequence of non-zero vectors $(x_n)$ in $E$, $T|_{[x_n]}$ is
  strictly singular.
\end{enumerate}
\end{Prop}
\begin{proof}
(1)$\Rightarrow$(2) Suppose that $M$ is a non-dispersed closed subspace of $E$ such that
$T|_M$ is an isomorphism.
Then there exists a normalized disjoint sequence $(x_k)$ in $E$ with
$\lim_{k\to\infty}\dist(x_k,M)=0$.
As in the proof of Lemma \ref{lemma:perturba}, we take a bounded sequence $(x^*_n)$ in $E^*$
such that $x^*_i(x_j)=\delta_{i,j}$ and, passing to a subsequence of $(x_k)$, we can find
a sequence $(m_k)$ in $M$ so that $\sum_{k=1}^\infty \|x^*_k\|\cdot\|x_k-m_k\|<1/2$.
Thus $Kx=\sum_{k=1}^\infty x^*_k(x) (x_k-m_k)$ defines $K \in\Lc(E)$ with $\|K\|<1/2$.
Hence $I-K$ is an isomorphism and $(I-K)x_k =m_k$.

Denoting $S=T(I-K)$, since $(I-K)([x_k]) =[m_k]$ and $T|_{[m_k]}$ is an isomorphism,
$S|_{[x_k]}$ is also an isomorphism, thus $T|_{[x_k]}=S|_{[x_k]}+TK|_{[x_k]}$ is upper
semi-Fredholm because $K$ is compact.
Since $N(T|_{[x_k]})$ is finite dimensional, Lemma \ref{lemma:fin} implies that
$T|_{[x_k : k\geq k_0]}$ is a isomorphism for some $k_0\in\N$, hence $T$ is not disjointly
strictly singular.

(2)$\Rightarrow$(3) Suppose that $(x_n)$ is a disjoint sequence of non-zero vectors in $E$ and
$T|_{[x_n]}$ is not strictly singular.
Then there exists an infinite dimensional closed subspace $M$ of $[x_n]$ such that $T|_M$ is
an isomorphism, and $M$ is not dispersed because there exists a normalized block basis $(y_k)$
of $(x_n)$ such that $\dist(y_k,M)<2^{-k}$.

(3)$\Rightarrow$(1) Strictly singular operators are never isomorphisms on infinite dimensional
subspaces.
\end{proof}

\begin{Def}
We say that $T\in\Lc(E,Y)$ is \emph{disjointly non-singular} if there is no disjoint sequence
of non-zero vectors $(x_n)$ in $E$ such that $T|_{[x_n]}$ is strictly singular.
\end{Def}

We denote $\dns(E,Y):=\{T\in\Lc(E,Y) : T \textrm{ is disjointly non-singular}\}$.

\begin{Thm}\label{thm:main}
For $T\in\Lc(E,Y)$, the following assertions are equivalent:
\begin{enumerate}
  \item $T$ is disjointly non-singular.
  \item There is no disjoint sequence of non-zero vectors $(x_n)$ in $E$ such that $T|_{[x_n]}$
  is compact.
  \item For every disjoint sequence of non-zero vectors $(x_n)$ in $E$, $T|_{[x_n]}$ is upper
  semi-Fredholm.
  \item For every normalized disjoint sequence $(x_n)$ in $E$, $\liminf_{n\to\infty}\|Tx_n\|>0$.
\end{enumerate}
\end{Thm}
\begin{proof}
(1)$\Rightarrow$(2) is immediate, and for (3)$\Rightarrow$(4) observe that, if $T|_{[x_n]}$
is upper semi-Fredholm, then Lemma \ref{lemma:fin} implies that $T|_{[x_n: n>n_0]}$ is an
isomorphism for some $n_0$.

(2)$\Rightarrow$(3) If $(x_n)$ is a disjoint sequence of non-zero vectors in $E$ and
$T|_{[x_n]}$ is not upper semi-Fredholm, then there exists an infinite dimensional closed
subspace $M$ of $[x_n]$ such that $T|_M$ is compact \cite[Theorem 7.16]{Aiena:04}.
Taking a normalized block basis $(y_k)$ of $(x_n)$ such that $\dist(y_k,M)<2^{-k}$, the
argument of Lemma \ref{lemma:perturba} allows us to show that $T|_{[y_k]}$ is compact.

(4)$\Rightarrow$(1) If $(x_n)$ is a disjoint sequence of non-zero vectors in $E$ such that
$T|_{[x_n]}$ is strictly singular, then we can construct a normalized block basis $(y_k)$ of
$(x_n)$ such that $\lim_{k\to\infty}\|Ty_k\|=0$.
\end{proof}

It was proved in \cite[Theorem 2]{GM:97} that $T\in\Lc(L_1,Y)$ is tauberian if and only
if it satisfies (4) in Theorem \ref{thm:main}.

\begin{Cor}\label{cor:pert}
Let $T,K\in\Lc(E,Y)$.
\begin{enumerate}
  \item If $T\in \dns$ and $K\in \dss$ then $T+K\in \dns$.
  \item  If $T\in \dns$ then $N(T)$ is dispersed.
\end{enumerate}
\end{Cor}
\begin{proof}
(1) is a consequence of Proposition \ref{prop:DSS-SS}, Theorem \ref{thm:main} and the fact
that the class of upper semi-Fredholm operators is stable under the addition of strictly
singular operators \cite[Theorem 7.46]{Aiena:04}.

(2) If $T\in \dns$ and $(x_n)$ is a disjoint sequence of non-zero vectors in $E$ then
$T|_{[x_n]}$ is upper semi-Fredholm.
Hence $[x_n]\cap N(T)$ is finite dimensional and $[x_n]+N(T)$ is closed, and the result
follows from Proposition \ref{prop:dispersed}.
\end{proof}

Next we give a perturbative characterization of disjointly non-singular operators.
Similar characterizations can be found in \cite[Theorem 7.16]{Aiena:04} for upper semi-Fredholm
operators, and in \cite{GO:90} for tauberian operators.

\begin{Thm}[Perturbative characterization]\label{thm:pert-char}
An operator $T\in\Lc(E,Y)$ is disjointly non-singular if and only if $N(T+K)$ is dispersed
for every compact operator $K\in\Lc(E,Y)$.
\end{Thm}
\begin{proof}
The direct implication is a consequence of Corollary \ref{cor:pert}.

For the converse, suppose that $T$ is not disjointly non-singular.
By Theorem \ref{thm:main} there exists a normalized disjoint sequence $(x_n)$ in $E$ such
that $\|Tx_n\|<2^{-n}$.
Since $(x_n)$ is a normalized basic sequence, there exists a bounded sequence $(x^*_n)$
in $E^*$ such that $x^*_i(x_j)=\delta_{ij}$ for $i,j\in\N$.

The operator $K:E\to Y$ defined by $Kx =-\sum_{n=1}^\infty x^*_n(x) Tx_n$ is compact, and
$[x_n]$ is contained in $N(T+K)$, hence $N(T+K)$ is not dispersed.
\end{proof}

A Banach lattice is called \emph{atomic} when its Banach lattice structure is induced by a
(countable or uncountable) $1$-unconditional basis.
For example, the subspace generated by a disjoint sequence of non-zero vectors of $E$ is an
atomic sublattice of $E$.

\begin{Cor}
If $E$ is atomic, then every operator in $\dns(E,Y)$ is upper semi-Fredholm
and every operator in $\dss(E,Y)$ is strictly singular.
\end{Cor}
\begin{proof}
We prove the countable case, since the proof of the general case is similar.
Let $(x_n)$ be a $1$-unconditional basis of $E$ inducing its Banach lattice structure.
A standard block-basis argument shows that for every closed infinite dimensional subspace $M$
of $E$ there is a normalized block basis $(y_k)$ such that $\lim_{k\to\infty}\dist(y_k,M)=0$.
Since $(y_k)$ is disjoint, we get that dispersed subspaces of $E$ are finite dimensional.
Thus the perturbative characterizations of the upper semi-Fredholm operators and $\dns(E,Y)$
(\cite[Theorem 7.16]{Aiena:04} and Theorem \ref{thm:pert-char}) imply the first part, and
Proposition \ref{prop:DSS-SS} implies the second part.
\end{proof}

For operators with closed range there is a simpler characterization of $T\in \dns$.

\begin{Prop}\label{prop:closed-range}
Let $T\in\Lc(E,Y)$ with closed range.
Then $T$ is disjointly non-singular if and only if $N(T)$ is dispersed.
\end{Prop}
\begin{proof}
The direct implication is proved in Corollary \ref{cor:pert}.
For the converse one, suppose that $T$ is not disjointly non-singular.
By Theorem \ref{thm:main} there exists a normalized disjoint sequence $(x_n)$ in $E$
such that $\|Tx_n\|<2^{-n}$.
Since $R(T)$ is closed there exists a constant $C>0$ so that $\dist(x_n,N(T))<C\cdot 2^{-n}$,
hence $N(T)$ is not dispersed.
\end{proof}

In the case $\dns(L_1,Y)\neq\emptyset$, it was proved in \cite[Proposition 14]{GM:97} that
$K\in\Lc(L_1,Y)$ is $\ell_1$-singular if and only if $T+K\in\dns$ for every $T\in\dns(L_1,Y)$.
This means that $\dss(L_1,Y)$ is the perturbation class of $\dns(L_1,Y)$.
Moreover $\dns(L_1,Y)$ is an open subset of $\Lc(L_1,Y)$.

\begin{Question}\label{Q1}
Suppose that $\dns(E,Y)$ is non-empty.
\begin{itemize}
  \item[(a)] Is $\dss(E,Y)$ the perturbation class of $\dns(E,Y)$?
  \item[(b)] Is $\dns(E,Y)$ an open subset of $\Lc(E,Y)$?
\end{itemize}
\end{Question}

In general, $T\in\Lc(X,Y)$ tauberian does not imply $T^{**}\in\Lc(X^{**},Y^{**})$ tauberian
\cite{AG:91}, but the implication is valid for $X=L_1$ \cite[Corollary 9]{GM:97}.

Also each operator $T\in\Lc(X,Y)$ induces an operator $T^{co}\in\Lc(X^{**}/X,Y^{**}/Y)$,
defined by $T^{co}(x^{**}+X)=T^{**}x^{**}+Y$, which is called the residuum operator.
An operator $T\in\Lc(L_1,Y)$ is tauberian if an only if $T^{co}$ is an isomorphism
\cite[Proposition 11]{GM:97}.

\begin{Question}
Suppose that $E$ is non-reflexive and $T\in \dns(E,Y)$.
\begin{itemize}
  \item[(a)] Is $T^{**}\in \dns$?
  \item[(b)] Is $T^{co}$ an isomorphism?
\end{itemize}
\end{Question}

For information on the residuum operator $T^{co}$ we refer to \cite[Section 3.1]{GM:10}
and \cite{GSK:95}.

\section{Disjointly non-singular operators on $L_p$}

Here we study the disjointly non-singular operators on $L_p$ ($1<p<\infty$).
The case $p=1$ was studied in \cite{GM:97}.

For $T\in\Lc(E,Y)$ we consider the following quantity:
$$
\beta(T) :=\inf\big\{\liminf_{n\to\infty}\|Tx_n\|: (x_n)
\textrm{ normalized disjoint sequence in }E\big\},
$$
and for a measurable function $f$ on $(0,1)$ we denote $D(f)=\{t\in (0,1) : f(t)\neq 0\}$.

\begin{Prop}
An operator $T\in\Lc(L_p,Y)$ ($1< p<\infty$) is disjointly non-singular if and only if
$\beta(T)>0$.
\end{Prop}
\begin{proof}
Clearly $\beta(T)>0$ implies (4) in Theorem \ref{thm:main}.

Conversely, suppose that $\beta(T)=0$.
Then for every $k\in\N$ there exists a normalized disjoint sequence $(f^k_n)_{n\in\N}$ in $L_p$
such that $\|Tf^k_n\|<1/k$ for every $n\in\N.$

We denote $D(k,n)=\{t\in (0,1): f^k_n(t)\neq 0\}$, and we take $n_1=1$ and $g_1= f^1_1$.
Since $\lim_{n\to\infty} \mu\left(D(k,n)\right)=0$, we can select $n_2>1$ such that
$\|g_1\chi_{D(2,n_2)}\|^p_p<2^{-2}$, and take $g_2= f^2_{n_2}$.
Proceeding in this way we obtain $1=n_1<n_2<\cdots$ in $\N$ so that $g_k= f^k_{n_k}$ satisfies
$\|g_i\chi_{D(g_k)}\|^p_p<2^{-2k}$ for $i<k$.
We denote $A_k=D(g_k)\setminus\cup_{n>k}D(g_n)$.

Taking $f_k= g_k\chi_{A_k}$ we obtain a disjoint sequence $(f_k)$ satisfying
$$
1\geq \|f_k\|^p_p \geq 1-\sum_{n>k}\|g_k\chi_{D(g_n)}\|_p^p\geq 1-\sum_{n>k}2^{-2k}\geq 1/2
$$
and
$$
\|Tf_k\| \leq \|Tg_k\| + \|T\|\sum_{n>k}\|g_k\chi_{D(g_n)}\|_p\leq 1/k +2^{-kp}\|T\|.
$$
Thus $T$ fails (4) in Theorem \ref{thm:main}.
\end{proof}

\begin{Cor}
$\dns(L_p,Y)$ is an open subset of $\Lc(L_p,Y)$.
\end{Cor}
\begin{proof}
Observe that $\beta(T+A)>\beta(T)-\|A\|$.
\end{proof}

\begin{Question}
Is $\beta(T)>0$ for every $T\in\dns(E,Y)$?
\end{Question}

We say that a closed subspace $M$ of $L_p$ ($1\leq p\leq \infty$) is \emph{strongly embedded}
if, in $M$, convergence in measure is equivalent to convergence in $\|\cdot\|_p$-norm.

The following result is proved in \cite[Proposition 6.4.5]{A-Kalton:06} for $p<\infty$,
but the proof is also valid for $p=\infty$.

\begin{Prop} \label{prop:strongly-embed}
A closed subspace $M$ of $L_p$ ($1\leq p\leq \infty$) is strongly embedded if and only if for
each (equivalently, for some) $0<q<p$, $\|\cdot\|_q$ and $\|\cdot\|_p$ are equivalent on $M$.
\end{Prop}

Next result describes the strongly embedded subspaces in $L_p$ for $p<\infty$, showing
that the cases $p>2$, $p=2$ and $p<2$ are remarkably different.

\begin{Prop}\label{prop:strongly-Lp}
Let $1\leq p<\infty$.
\begin{enumerate}
  \item For $p\neq 2$, a closed subspace of $L_p$ is strongly embedded if and only if it
  contains no subspace isomorphic to $\ell_p$.
  \item For $p>2$, a closed infinite dimensional subspace of $L_p$ is strongly embedded if
  and only it is isomorphic to $\ell_2$.
  \item For $p<2$, the set of strongly embedded subspaces of $L_p$ include isomorphic copies of
  $L_q$ for every $q\in (p,2]$.
  \item There exists an orthogonal decomposition $L_2=M\oplus M^\perp$ with both $M$ and
  $M^\perp$ strongly embedded in $L_2$.
  \item For $p\neq 2$, we cannot write $L_p$ as the direct sum of two strongly embedded
subspaces.
\end{enumerate}
\end{Prop}
\begin{proof}
We refer to \cite[Theorems 6.4.8 and 7.2.6]{A-Kalton:06} for (1), \cite[Theorem 6.4.8]{A-Kalton:06}
for (2), and \cite[Corollary 2.f.5]{LT-II:79} for (3).
Moreover, (4) is \cite[Theorem 8.22]{Pisier:CBMS}, and (5) follows from the fact that containing no
copy of $\ell_p$ is stable under direct sums.
\end{proof}

For additional information on strongly embedded subspaces of $L_p$ we refer to \cite{Astashkin:14},
where they are called $\Lambda(p)$-spaces.

Next result is a special case of the Kadec-Pe\l czy\'{n}ski dichotomy as stated in
\cite[Proposition 1.c.8]{LT-II:79}.
We give an alternative proof.

\begin{Prop}\label{prop:dispersed-Lp}
For $1\leq p<\infty$, a closed subspace $M$ of $L_p$ is strongly embedded if and only
if it is dispersed.
\end{Prop}
\begin{proof}
Suppose that $M$ is strongly embedded, and let $(f_n)$ be a disjoint sequence of non-zero vectors
in $L_p$.
For $p\neq 2$, since $M$ contains no copy of $\ell_p$ (Proposition \ref{prop:strongly-Lp}) and
$[f_n]$ is isomorphic to $\ell_p$, $M$ and $[f_n]$ are totally incomparable, hence $M\cap [f_n]$
is finite dimensional and $M+[f_n]$ is closed (see \cite{GO:86}).
Thus $M$ is dispersed by Proposition \ref{prop:dispersed}.

For $p=2$, it follows from Proposition \ref{prop:strongly-embed} that $\|\cdot\|_2$ and
$\|\cdot\|_1$ are equivalent on $M$.
But, since every operator in $\Lc(\ell_2,\ell_1)$ is compact \cite[Theorem 2.1.4]{A-Kalton:06},
the norms $\|\cdot\|_2$ and $\|\cdot\|_1$ are equivalent in no infinite dimensional subspace of
$[f_n]$, hence $M\cap [f_n]$ is finite dimensional.
Moreover, the argument in the proof of Lemma \ref{lemma:perturba} allows us to show that $M+[f_n]$
is closed.
\medskip

Conversely, if $M$ is not strongly embedded, arguing as in the proof of \cite[Theorem 6.4.7]{A-Kalton:06},
we get a normalized sequence $(g_n)$ in $M$ and a sequence of disjoint sets $(A_n)$ such that
$\|g_n-g_n\chi_{A_n}\|_p\to 0$.
Then $(g_n)$ is almost disjoint, hence $M$ is not dispersed.
\end{proof}

\begin{Rem}
\emph{In $L_\infty$, there are dispersed subspaces which are not strongly embedded.}

\emph{Indeed, strongly embedded subspaces of $L_\infty$ are finite dimensional because they are
reflexive and, for $1\leq p<\infty$, the natural map from $L_\infty$ to $L_p$ takes weakly
convergent sequences into convergent sequences \cite[Theorem 5.4.5]{A-Kalton:06}.
See \cite[Theorem 5.2]{Rudin:91} for an alternative argument.
However, each closed subspace of $L_\infty$ containing no copy of $c_0$ is dispersed because
it is totally incomparable with the subspace generated by a disjoint sequence of non-zero vectors.
See the proof of Proposition \ref{prop:dispersed-Lp}.}
\end{Rem}

The following result is essentially known, but we give a proof of it for convenience.

\begin{Prop}
Let $1\leq p<\infty$, $p\neq 2$.
An operator $T\in\Lc(L_p,Y)$ is disjointly strictly singular if and only if it is
$\ell_p$-singular.
\end{Prop}
\begin{proof}
By Proposition \ref{prop:DSS-SS}, $T$ is disjointly strictly singular if and only given
a closed subspace $M$ of $L_p$, $T|_M$ isomorphism implies $M$ dispersed.
Moreover, Propositions \ref{prop:strongly-Lp} and \ref{prop:dispersed-Lp} show that $M$ is
dispersed in $L_p$ ($1\leq p<\infty$, $p\neq 2$) if and only if $M$ contains no copy of
$\ell_p$.
These two facts imply the result.
\end{proof}

Note that $T\in\Lc(L_2,Y)$ is $\ell_2$-singular if and only it is strictly singular.

\begin{Prop}\label{prop:Lp-neq2}
Let $1< p<\infty$, $p\neq 2$.
For $T\in\Lc(L_p,Y)$, the following assertions are equivalent:
\begin{enumerate}
  \item $T$ is disjointly non-singular.
  \item $T|_M$ is upper semi-Fredholm for every subspace $M$ of $L_p$ isomorphic to $\ell_p$.
  \item $T|_M$ is compact for no subspace $M$ of $L_p$ isomorphic to $\ell_p$.
\end{enumerate}
\end{Prop}
\begin{proof}
(2)$\Leftrightarrow$(3) The direct implication is trivial.
Conversely, if there is a subspace $M$ of $L_p$ isomorphic to $\ell_p$ such that
$T|_M$ is not upper semi-Fredholm, then $M$ contains a closed infinite dimensional subspace
$N$ such that $T|_N$ is compact, and $N$ contains a subspace isomorphic to $\ell_p$.

(1)$\Rightarrow$(3) Suppose that there is a subspace $M$ of $L_p$ isomorphic to $\ell_p$ such
that $T|_M$ is compact.
By (1) in Proposition \ref{prop:strongly-Lp}, $M$ is not strongly embedded; hence there is a
normalized disjoint sequence $(f_n)$ in $L_p$ with $\lim_{n\to\infty}\dist(f_n,M)=0$ by
Proposition \ref{prop:dispersed-Lp}.
Passing to a subsequence, Lemma \ref{lemma:perturba} allows us to show that $T|_{[f_n]}$ is
compact, hence $T$ is not disjointly non-singular.

(2)$\Rightarrow$(1) It is a consequence of Theorem \ref{thm:main}, because the closed subspace
generated by a disjoint sequence of non-zero vectors in $L_p$ is isomorphic to $\ell_p$.
\end{proof}

For $p=2$, each of the assertions (2) and (3) in Proposition \ref{prop:Lp-neq2} is
equivalent to $T$ being upper semi-Fredholm.
\medskip

The following result shows some differences between the properties of $\dns(L_2)$ and those of
$\dns(L_p)$ for $p\neq 2$.

\begin{Prop}
Let $1<p<\infty$, $p\neq 2$.
\begin{enumerate}
  \item $S, T\in \dns(L_p)$ implies $ST\in \dns(L_p)$.
  \item There exist $S, T\in \dns(L_2)$ such that $ST=0$.
\end{enumerate}
\end{Prop}
\begin{proof}
(1) It follows from Proposition \ref{prop:Lp-neq2}: $T\in \dns(L_p)$ if and only if
$T|_M$ is upper semi-Fredholm for every subspace $M$ isomorphic to $\ell_p$.

(2) We consider the decomposition $L_2=M\oplus M^\perp$ with both $M$ and $M^\perp$ strongly
embedded given in Proposition \ref{prop:strongly-Lp}.
Let $S$ denote the orthogonal projection on $L_2$ onto $M$ and let $T=I-S $.
Since both $S$ and $T$ have closed range and dispersed kernel, $S, T\in \dns(L_2)$ by
Proposition \ref{prop:closed-range}.
\end{proof}

Next we show that the dispersed subspaces of $L_p$ form an open subset of the set of all
closed subspaces with respect to the gap metric \cite[Chapter IV]{Kato:80}.

Let $M$ and $N$ be closed subspaces of $X$, and let us denote $S_M=\{m\in M : \|m\|=1\}$.
The \emph{gap between $M$ and $N$} is defined by
$g(M,N)=\max \big\{\delta(M,N), \delta(N,M)\big\}$,
where $\delta(M,N)=\sup_{m\in S_M}\dist(m,N)$.

\begin{Prop}\label{prop:stability}
Given a dispersed subspace $M$ of $L_p$ ($1\leq p<\infty$), there exists $\e>0$ such that
if $M_1$ is a closed subspace of $L_p$ and $\delta(M_1,M)<\e$, then $M_1$ is dispersed.
\end{Prop}
\begin{proof}
By Proposition \ref{prop:closed-range} and \cite{GM:97} (case $p=1$), the quotient map
$Q:L_p\to L_p/M$ is in $\dns$, hence $\beta(Q)>0$.
Thus, for every normalized disjoint sequence $(x_n)$ in $L_p$,
$$
\liminf_{n\to\infty}\|Q x_n\|=\liminf_{n\to\infty}\dist(x_n,M)\geq \beta(Q)>0.
$$
By \cite[Lemma IV.4.2]{Kato:80},
$\big(1+\delta(M_1,M)\big)\dist(x_n,M_1)\geq \dist(x_n,M)- \delta(M_1,M)$.
Thus we can take $\e=\beta(Q)/2$.
\end{proof}

As a consequence, we obtain a stability result for the decompositions of $L_2$ as a direct
sum of strongly embedded subspaces, like the one given in Proposition \ref{prop:strongly-Lp}.

\begin{Cor}\label{cor:stability}
Let $M$ and $N$ be strongly embedded subspaces of $L_2$ with $L_2=M\oplus N$.
Then there exists $\e>0$ such that if $M_1$ and $N_1$ are closed subspaces of $L_2$,
$g(M_1,M)<\e$ and $g(N_1,N)<\e$, then $M_1$ and $N_1$ are strongly embedded and
$L_2=M_1\oplus N_1$.
\end{Cor}
\begin{proof}
Applying \cite[Theorem IV.4.24]{Kato:80} twice, we can find $\e>0$ such that $g(M_1,M)<\e$
and $g(N_1,N)<\e$ imply $L_2=M_1\oplus N_1$.
And, by Proposition \ref{prop:stability}, we can choose $\e$ so that additionally $M_1$ and
$N_1$ are strongly embedded.
\end{proof}

There are some other differences between $\dns(L_2)$ and $\dns(L_p)$ for $p\neq 2$.

\begin{Prop}\label{prop:beta}
Let $1< p<2$ and let $M$ be a dispersed subspace of $L_p$.
Then the quotient map $Q_M:L_p\to L_p/M$ satisfies $\beta(Q_M)=1$.
\end{Prop}
\begin{proof}
We showed in Proposition \ref{prop:dispersed-Lp} that, for $1\leq p<\infty$, a closed
subspace $L_p$ is strongly embedded if and only if it is dispersed.
Moreover, for $1<p<2$, it was proved in \cite[Theorem 2]{Astashkin:14} that the unit ball $B_M$
of a strongly embedded subspace of $L_p$ is equi-integrable in $L_p$; i.e., that for each $\e>0$
there exists $\delta >0$ such that if $f\in B_M$ and $A$ is a measurable subset of $(0,1)$ with
$\mu(A)<\delta$ then $\|f\cdot\chi_A\|_p<\e$.

If $f\in L_p$ and $\|f\|_p=1$ then $\|Q_M f\| = \dist(f,M)= \inf_{g\in 3 B_M} \|f-g\|_p$.
Thus, given a normalized disjoint sequence $(f_n)$ in $L_p$ and $g\in 3 B_M$, and denoting
$A_n=D(f_n)$, since $\lim_{n\to\infty}\mu(A_n)=0$, there exists $n_\e$ so that
$$
\|f_n-g\|_p \geq \|f_n-g\cdot\chi_{A_n}\|_p\geq 1-\e
$$
for $n\geq n_\e$; hence $\liminf_{n\to\infty}\|Q_M f_n\|=1$ and we get $\beta(Q_M)=1$.
\end{proof}

\begin{Rem}
\emph{Proposition \ref{prop:beta} fails for $p=2$.}

\emph{Indeed, if we consider the decomposition $L_2=M\oplus M^\perp$ with both $M$ and $M^\perp$
strongly embedded given in Proposition \ref{prop:strongly-Lp}, then $\beta(Q_M)>0$ and
$\beta(Q_{M^\perp})>0$.
Moreover, for $f\in L_2$ we have $\|f\|_2^2= \dist(f,M)^2+\dist(f,M^\perp)^2$, which implies
$\beta(Q_M)^2 +\beta(Q_{M^\perp})^2\leq 1$. Hence $0<\beta(Q_M), \beta(Q_{M^\perp})<1$.}
\end{Rem}

\begin{Question}
Let $M$ be a dispersed subspace of $L_p$ ($2< p<\infty$).
Is $\beta(Q_M)=1$?

For $2< p<\infty$ there are strongly embedded subspaces of $L_p$ whose unit ball is not
equi-integrable in $L_p$ \cite{Astashkin:14}.
So the argument in the proof of Proposition \ref{prop:beta} is not valid.
\end{Question}

Next we give further characterizations of $\dns(L_p,Y)$ and some consequences.

\begin{Thm}\label{thm:Lp-D}
For $T\in\Lc(L_p,Y)$, $1< p<\infty$, the following assertions are equivalent:
\begin{enumerate}
  \item $T$ is disjointly non-singular.
  \item For every normalized sequence $(f_n)$ in $L_p$ with $\lim_{n\to\infty}\mu\left(D(f_n)\right)=0$,
  we have $\liminf_{n\to\infty}\|Tf_n\|>0$.
  \item There exists $r>0$ such that $f\in L_p$, $\|f\|_p=1$ and $\mu\left(D(f)\right)\leq r$
  imply $\|Tf\|\geq r$.
\end{enumerate}
\end{Thm}
\begin{proof}
(3)$\Rightarrow$(2) is clear and (2)$\Rightarrow$(1) follows from Theorem \ref{thm:main}.

For the remaining implication, suppose that (3) fails.
Then we can find a normalized sequence $(f_n)$ in $L_p$ such that $\mu\left(D(f_n)\right)< 1/n$
and $\|Tf_n\|<1/n$.
Passing to a subsequence we can find a normalized disjoint sequence $(g_n)$ with
$\lim_{n\to\infty}\|f_n-g_n\|=0$ and $\lim_{n\to\infty}\|Tg_n\|=0$, hence (1) fails.
\end{proof}

\begin{Cor}\label{cor:cont-Lp}
Let $1< p<\infty$. If $\dns(L_p,Y)$ is non-empty, then $Y$ contains a subspace
isomorphic to $L_p$.
\end{Cor}
\begin{proof}
Let $T\in \dns(L_p,Y)$ and let $r$ be as in Theorem \ref{thm:Lp-D}.
Note that we can assume that $0<r<1$, hence the restriction $T|_{L_p(0,r)}$ is an isomorphism.
\end{proof}

Let us see that the perturbation classes problem (Question \ref{Q1}) has a positive answer
for $E=L_p$ ($1<p<\infty$).
The case $p=1$ was proved in \cite{GM:97}.

\begin{Thm}\label{thm:pert-class-Lp}
Let $1< p<\infty$ and suppose that $\dns(L_p,Y)$ is non-empty.
An operator $K\in \Lc(L_p,Y)$ is in $\dss$ if and only if $T+K\in \dns$ for every
$T\in \dns(L_p,Y)$.
\end{Thm}
\begin{proof}
The direct implication is a consequence of Corollary \ref{cor:pert}.

The proof of the converse implication is similar to that of \cite[Proposition 14]{GM:97} for
the case $p=1$.
We suppose that $K\notin \dss$, and we will construct $T\in \dns(L_p,Y)$ such that
$T+K\notin \dns(L_p,Y)$.

Since $K\notin \dss$, there exists a disjoint sequence of non-zero vectors $(f_n)$ in $L_p$ such
that $K|_{[f_n]}$ is an isomorphism.
By \cite[Proposition 6.4.1]{A-Kalton:06}, $[f_n]$ is a complemented subspace of $L_p$ isomorphic
to $\ell_p$ and $N:=K([f_n])$ is a subspace of $Y$ isomorphic to $\ell_p$.
By Corollary \ref{cor:cont-Lp}, $Y$ contains a subspace $L$ isomorphic to $L_p$.

Let $M$ be a closed complement of $[f_n]$ in $L_p$, and let $U:M\to L$ be an isomorphic embedding.
Considering the relative positions of the subspaces $L$ and $N$ inside $Y$, we have three cases:

(a) \emph{$L\cap N$ finite dimensional and $L+N$ closed.}
By Lemma \ref{lemma:fin} we can assume that $L\cap N=\{0\}$.
Then $T:L_p=M\oplus [f_n]\to Y$ defined by $T|_M = U$ and $T|_{[f_n]}= -K|_{[f_n]}$ is an
isomorphic embedding, hence $T\in\dns(L_p,Y)$, but $T+K\notin \dns$ because $(f_n)\subset N(T+K)$.

(b) \emph{$L\cap N$ infinite dimensional.}
In this case $L\cap N$ contains a subspace $N_1$ isomorphic to $\ell_p$ and complemented in $L$
\cite[Theorems 6.4.8 and 7.2.6]{A-Kalton:06}, hence $M_1= (K|_M)^{-1}(N_1)$ is complemented in
$L_p$ and isomorphic to $\ell_p$.
Since $L_p$ is primary \cite[Theorem 2.d.11]{LT-II:79}, the complements of $N_1$ and $M_1$ in $L$
and $L_p$ are isomorphic to $L_p$, so we can construct $T$ as in the previous case.

(c) \emph{$L\cap N$ finite dimensional and $L+N$ non-closed.}
An argument in \cite[Proof of Theorem 4.3.5]{GM:10} provides a compact operator
$K_1\in \Lc(L_p,Y)$ such that $(K+K_1)|_{[f_n]}$ is an isomorphism and $L\cap (K+K_1)([f_n])$
is infinite dimensional.
The argument in case (b) gives $T\in \dns(L_p,Y)$ such that $T+K+K_1\notin \dns$, hence
$T+K\notin \dns$.
\end{proof}

Each upper semi-Fredholm operator is tauberian, but no other examples of tauberian operators
in $\Lc(L_1)$ were known (see \cite[Section 4.1]{GM:10}) until, using probabilistic arguments,
examples with infinite dimensional kernel or non-closed range were constructed in \cite{JNST:14}.
Next we show that it is much easier to find similar examples in $\Lc(L_p)$, $1<p<\infty$.

\begin{Example}
\emph{For $1<p<\infty$, there exists a projection $P\in\Lc(L_p)$ onto the closed subspace
generated by the Rademacher functions.
See the proof of \cite[Proposition 6.4.2]{A-Kalton:06}.}

\emph{The operator $I-P:L_p\to L_p$ has closed range and infinite dimensional kernel.
In fact, $N(I-P)$ is the closed subspace generated by the Rademacher functions, which is
strongly embedded \cite[Proposition 6.4.5]{A-Kalton:06}.
Thus the operator $I-P$ is disjointly non-singular by Proposition \ref{prop:closed-range}.}

\emph{Also, it is not difficult to find a compact operator $K\in\Lc(L_p)$ such that the range
of $I-P+K$ is non-closed.
Note that $I-P+K$ is disjointly non-singular.}
\end{Example}

\end{document}